\documentclass[11pt,a4paper,reqno]{amsart}
\usepackage[T1]{fontenc}
\usepackage{lmodern}
\usepackage{mathrsfs}
\usepackage{amsmath,amssymb,amsfonts,amsthm}
\usepackage{latexsym,bm}
\usepackage{fancyhdr,graphicx}
\usepackage{indentfirst}
\usepackage{enumerate}
\usepackage{url}
\numberwithin{equation}{section}
\newtheorem{theorem}{Theorem}[section]
\newtheorem{proposition}[theorem]{Proposition}
\newtheorem{corollary}[theorem]{Corollary}
\newtheorem{lemma}[theorem]{Lemma}
\theoremstyle{definition}
\newtheorem{definition}[theorem]{Definition}

\makeatletter
\@ifundefined{subjclassname@2020}{\@namedef{subjclassname@2020}{\textup{2020} Mathematics Subject Classification}}{}
\makeatother
\subjclass[2020]{34C07, 34C08, 37G15}
\keywords{Quadratic center, period annulus, cyclicity, Abelian integral, Picard--Fuchs equation, Chebyshev system}
\begin{document}
\title[Cyclicity of quadratic codimension-four centers]{Cyclicity of period annuli of quadratic codimension four centers}
\author{Qi Meng}
\author{Yulin Zhao}
\address{School of Mathematics (Zhuhai), Sun Yat-sen University, Zhuhai Campus, Zhuhai, 519082, China}
\email{mengq25@mail2.sysu.edu.cn (Q. Meng); mcszyl@mail.sysu.edu.cn (Y. Zhao)}
\date{}
\begin{abstract}
Quadratic systems with a center fall into four classes:
Hamiltonian $Q_3^H$, reversible $Q_3^R$, generalized
Lotka--Volterra $Q_3^{LV}$, and codimension-four $Q_4$.
Before this work, the cyclicity of period annuli under
quadratic perturbations had been completely determined for
$Q_3^H$, whereas only certain special cases had been resolved
for $Q_3^R$, $Q_3^{LV}$ and $Q_4$.
In this paper, we establish that the period annulus of every
center in $Q_4$ has cyclicity three under quadratic
perturbations, thereby proving \.{Z}o{\l}{\k{a}}dek's 1994
conjecture for this class.
The proof combines Picard--Fuchs equations with a Chebyshev
criterion. Its key step is a strict Wronskian inequality for
two ratios of elliptic integrals on their entire interval of
definition. We establish this inequality by constructing a
comparison function and excluding a first zero of the
Wronskian. The Chebyshev property avoids the two
additional zeros allowed by the second-order differential
reduction used in earlier proofs.
\end{abstract}
\maketitle

\section{Introduction}\label{sec:intro}

The second part of Hilbert's sixteenth problem asks for a bound,
depending only on the degree, for the number of limit cycles of a
planar polynomial vector field, and for their possible configurations.
Even for quadratic systems, the determination of the maximal number
remains unresolved. A natural approach is to study small perturbations of
integrable systems.

In this paper we study the number of limit cycles that can
bifurcate from period annuli of quadratic integrable systems
under quadratic perturbations.
For a prescribed class of perturbations, the \emph{cyclicity} of
a compact invariant set $K$ of the unperturbed vector field is
the supremum of the numbers of limit cycles that can occur
simultaneously in arbitrarily small neighborhoods of $K$ under
arbitrarily small perturbations. It is infinite if these numbers
have no finite upper bound.
Such compact configurations include singular points, periodic
orbits and polycycles. For a period annulus $\Pi$, we define
\[
\operatorname{Cycl}(\Pi)
 =\sup_{K\Subset\Pi}\operatorname{Cycl}(K),
\]
where $K$ runs over closed subannuli bounded by periodic orbits.
Thus the cyclicity of the center and that of the outer
singular boundary are considered separately.

Let $\Pi$ be a period annulus of an integrable system, and let
$H$ be a first integral with $\nabla H\ne0$ on $\Pi$. Write
\[
\Pi=\bigcup_{h_c<h<h_s}\gamma_h,\qquad
\gamma_h\subset\{(x,y):H(x,y)=h\},
\]
where $\gamma_h$ is a periodic orbit, and $h_c$ and $h_s$ are
the limiting values of $H(x,y)$ corresponding to the inner and outer
boundaries of $\Pi$, respectively.

Let $M(x,y)$ be a nonvanishing integrating factor of the
unperturbed system on $\Pi$. For a small parameter $\varepsilon$,
consider an analytic perturbation of the form
\[
\dot x=\frac{H_y(x,y)}{M(x,y)}
       +\varepsilon X_n(x,y,\varepsilon),\qquad
\dot y=-\frac{H_x(x,y)}{M(x,y)}
       +\varepsilon Y_n(x,y,\varepsilon),
\]
where $X_n(x,y,\varepsilon)$ and $Y_n(x,y,\varepsilon)$ are polynomials in $x,y$ of degree at
most $n$ with coefficients
analytic in $\varepsilon$.
Choose a transverse section meeting each $\gamma_h$ once and
use $h$ as its coordinate. The displacement
$d(h,\varepsilon)$ of the Poincar\'e return map on this section
has the expansion
\begin{equation*}
 d(h,\varepsilon)=\varepsilon M_1(h)
                  +\varepsilon^2M_2(h)+\cdots.
\end{equation*}
The coefficient $M_k(h)$ is called the $k$-th order Melnikov
function. If the displacement is not identically zero, let
$M_k(h)$ be its first nonzero Melnikov function. On compact
subintervals of $(h_c,h_s)$, the zeros of $M_k(h)$, counted with
multiplicity, bound the number of bifurcating limit cycles,
and a simple zero gives a hyperbolic limit cycle.
The first-order Melnikov function is
\[
M_1(h)=\oint_{\gamma_h}M(x,y)
 \bigl[Y_n(x,y,0)\,dx-X_n(x,y,0)\,dy\bigr],
\]
where $\gamma_h$ is oriented by the unperturbed flow.

In the polynomial Hamiltonian setting, $M(x,y)\equiv1$, and the
first-order Melnikov function for a polynomial perturbation
is an Abelian integral.
The weak, or infinitesimal, Hilbert sixteenth problem, proposed
by Arnold~\cite[Problems~1978-6 and~1979-16]{Arnold2004}, asks
for a uniform upper bound for the number of isolated zeros of such integrals in
terms of the degrees of the Hamiltonian and the perturbation.
For polynomial Hamiltonian systems of arbitrary degree,
Binyamini, Novikov and
Yakovenko~\cite{BNY} obtained an explicit double-exponential
upper bound. Grau, Ma\~nosas and Villadelprat~\cite{GMV}
developed Chebyshev criteria for Abelian integrals, and Liu and
Xiao~\cite{LiuXiao2013} subsequently obtained a monotonicity
criterion for ratios of two Abelian integrals associated with
potential Hamiltonians.

More generally, Huang, Liang and
Zhang~\cite{HuangLiangZhang2023} established
a general Chebyshev criterion for families of integrals,
with applications to first-order Melnikov functions of
integrable non-Hamiltonian systems.
For multi-parameter center families admitting an explicit
first integral or inverse integrating factor,
Itikawa, Oliveira and Torregrosa~\cite{ItikawaEtAl2022}
obtained a formula for the first-order Taylor expansion
of the first-order Melnikov function with respect to
the center parameters.
Further results can be found in these works and the references
therein.

Bounds for the number of zeros of Abelian integrals do not, by themselves,
determine the cyclicity of a period annulus. If $M_1(h)\equiv0$
but the displacement is not identically zero, the first nonzero
higher-order Melnikov function must be considered.
Gavrilov and Iliev~\cite{GavrilovIliev2023} used higher-order
Melnikov functions to bound the cyclicity of period annuli
in a generalized Rayleigh--Li\'enard oscillator.

We now restrict the discussion to quadratic perturbations of
quadratic systems with a center.
In the space of all quadratic systems, the systems with a
center form a union of four irreducible affine algebraic sets:
the Hamiltonian $Q_3^H$,
reversible $Q_3^R$, generalized Lotka--Volterra $Q_3^{LV}$ and
codimension four $Q_4$~\cite[Theorem~1]{Zoladek}.
Following~\cite[p.~69]{GI}, a center is called generic if it
belongs to only one of these sets.
To study the cyclicity of the period annuli surrounding
these quadratic centers, Iliev~\cite{Iliev1998} determined
the essential perturbations and computed the corresponding
bifurcation functions.
These perturbations form selected families that suffice to
determine the cyclicity of the corresponding period annuli
under arbitrary quadratic perturbations.

The cyclicity of period annuli of quadratic Hamiltonian
systems under quadratic perturbations has been completely
determined. The period annuli of generic Hamiltonian centers
have cyclicity two~\cite{HI,Gavrilov2001,LiZhang2002}.
Chen et al.~\cite{ChenEtAl} subsequently gave a unified proof
in the real domain of the corresponding sharp bound of two
zeros for Abelian integrals.

Hamiltonian centers in $Q_3^R$, i.e., those in
$Q_3^H\cap Q_3^R$, are called degenerate.
Iliev~\cite{IlievTriangle,Iliev1998}
obtained formulas for the higher-order Melnikov functions
used in subsequent studies of the cyclicity of period
annuli~\cite{GavrilovIliev2000,ZhaoLiangLu2000,ZhaoZhu2001}.
Chow, Li and Yi~\cite{ChowLiYi2002} resolved the remaining
elliptic-segment and saddle-loop cases.
The cyclicity is three for the Hamiltonian
triangle and two for all other degenerate cases.

For non-Hamiltonian centers in $Q_3^R$, quadratic perturbations
of families with two centers were studied
in~\cite{CollLiProhens2009,IlievLiYu2005,Liu2012}.
Ji and Liu~\cite{JiLiu2022} obtained sharp cyclicity bounds of
two or three, depending on the parameter, for a one-parameter
family of reversible centers with elliptic phase curves.
Recent results include an upper bound of two for a reversible
system with an exponential first integral~\cite{LiuLiLlibre2022}
and cyclicity two for a family with logarithmic first
integrals~\cite{CenEtAl2023}.
For the outer boundaries of period annuli, Mar\'in and
Villadelprat~\cite{MarinVilladelprat2025} studied the cyclicity
of hyperbolic hemicycles under quadratic perturbations.
Further results can be found in these works and the references
therein.

For generalized Lotka--Volterra centers in $Q_3^{LV}$,
\.{Z}o{\l}{\k{a}}dek~\cite[Theorem~5]{Zoladek} obtained the
sharp bound of two for the cyclicity of period annuli of
non-Hamiltonian centers with two real invariant lines.
Fran\c{c}oise and Gavrilov~\cite{FrancoiseGavrilov2022} used
iterated integrals of length two to study quadratic deformations
of the isochronous Lotka--Volterra double center; see also the
references therein.

For centers whose complexified orbits have genus one,
Gautier, Gavrilov and Iliev~\cite{GGI} developed a systematic
elliptic-integral approach. Picard--Fuchs equations, Riccati
equations and Chebyshev properties are central tools in these
investigations.

We now turn to $Q_4$. For its two non-generic
centers, the cyclicity of the period annulus is
three, by Iliev~\cite[Theorem~1.2]{Iliev1997} and
\cite[Theorem~5]{Iliev1998}.
For generic centers, Gavrilov and Iliev~\cite{GI} obtained
the upper bound eight for the cyclicity of the period annulus
by combining Picard--Fuchs
equations with the argument principle. In 2011,
Zhao~\cite{Zhao} improved the upper bound to five and established
the lower bound three. The conjectured sharp bound is three
\cite[Conjectures~A and~B, Fig.~7]{Zoladek}.
To the best of our knowledge, these bounds for generic
$Q_4$ period annuli have remained unimproved since 2011.

In this paper, we study the cyclicity of period annuli of
generic $Q_4$ centers under quadratic perturbations.
Combining our result for generic $Q_4$ centers with Iliev's
results for the non-generic cases gives the following theorem.

\begin{theorem}\label{thm:main}
The period annulus of every quadratic center in $Q_4$
has cyclicity three under
quadratic perturbations.
\end{theorem}

The lower bound three was already
established in~\cite[Theorem~1 and the proof of Corollary~9]{Zhao}
by constructing three limit cycles using the Taylor expansion
of the generating function at the center endpoint. Therefore,
to prove Theorem~\ref{thm:main}, it remains to establish the
upper bound three for generic $Q_4$ centers.
As a by-product of the upper-bound proof, we obtain
an alternative construction attaining the lower bound
three, as stated in the following theorem.

\begin{theorem}\label{thm:interpolation}
Let $\Pi$ be the period annulus of a quadratic system with a
generic center in $Q_4$. Given any three distinct periodic
orbits in $\Pi$, there exists an analytic one-parameter
quadratic perturbation of this system which, for all
sufficiently small nonzero values of its parameter, has three
hyperbolic limit cycles converging to the prescribed orbits
as the parameter tends to zero.
\end{theorem}

With Theorem~\ref{thm:main}, the cyclicity of period annuli
under quadratic perturbations has now been completely
determined for two of the four classes of quadratic centers,
$Q_3^H$ and $Q_4$. This settles the $Q_4$ case of
\.{Z}o{\l}{\k{a}}dek's 1994 conjecture~\cite{Zoladek}.
For $Q_3^R$ and $Q_3^{LV}$, the cyclicity problem remains
open in general.

The earlier proofs in~\cite{GI,Zhao} estimate the number of
zeros through a second-order differential reduction that
allows two additional zeros. We prove the extended complete
Chebyshev property of the four functions arising before this
reduction. The main new step is a strict Wronskian inequality
for two ratios of elliptic integrals throughout their interval
of definition. We establish this inequality using an explicit
comparison function and two polynomial sign identities to
exclude a first zero of the Wronskian. A differential identity and Rolle's theorem, together with
the vanishing of the generating function at the center endpoint,
then give the upper bound of three for the number of zeros of
each generating function that is not identically zero.

The paper is organized as follows. Section~\ref{sec:reduction}
introduces the generating function and outlines the reduction
to a Wronskian inequality.
Section~\ref{sec:prelim} collects the Picard--Fuchs and Riccati
equations, center expansions and preliminary results on
Chebyshev systems and zeros of analytic functions.
Section~\ref{sec:comparison} derives a second-order inhomogeneous
equation and constructs a comparison function.
Section~\ref{sec:positive} establishes positivity of the decisive
Wronskian. Section~\ref{sec:conclusion} completes the estimate for
the number of zeros of the generating function and the proofs
of Theorems~\ref{thm:main} and~\ref{thm:interpolation}.

\section{Reduction and outline of the proof}\label{sec:reduction}

In this section, we outline the proof in three steps.
First, we reduce the
cyclicity problem to estimating the number of zeros of the
generating function $I(h)$. Second, we reduce the latter
problem to estimating the number of zeros of linear
combinations of four auxiliary functions. Third, we state
the key Wronskian inequality and
explain how the Chebyshev criterion and Rolle's theorem
complete the argument.

\subsection{The generating function and its zeros}

In the complex coordinate $z=x+iy$, a generic $Q_4$ center
admits the normal form
\begin{equation*}
\dot z=-iz+4z^2+2|z|^2+(b+ic)\overline z^{\,2},
\qquad |b+ic|=2,\quad c\ne0.
\end{equation*}
The parameter
\[
\kappa=\frac4{2+b}>1
\]
covers these centers, up to the equivalences in~\cite{GI}.
The restriction $c\ne0$, equivalently $-2<b<2$, specifies the
generic case treated below. The two non-generic cases
$c=0$, $b=\pm2$, are covered by Iliev's results in the proof of
Theorem~\ref{thm:main}.

We recall the integral representation and retain the notation
of~\cite{GI,Zhao}.
For fixed $\kappa>1$,
put
\begin{equation*}
H(x,y)=\frac23(\kappa-1)x^3-(\kappa-1)x^2y
       +\frac\kappa3y^3-y.
\end{equation*}
The level curve $H(x,y)=h$ has an oval $\gamma_h$
surrounding $(1,1)$ for
\[
 h_c=-\frac23<h<h_s=-\frac{2}{3\sqrt\kappa}.
\]
Here $h_c=H(1,1)$ corresponds to the center, and $h_s$ to
the outer boundary of the period annulus.
Let $D_h$ be the domain enclosed by $\gamma_h$, and define
\[
I_{i,j}(h)=\iint_{D_h}x^iy^j\,dx\,dy.
\]
The polynomial $H(x,y)$ is obtained by the
integral transformations of~\cite{GI}; it does not identify the
original quadratic vector field with a cubic Hamiltonian vector
field.

Iliev's essential-perturbation theorem~\cite[Theorems~1
and~2(iii)]{Iliev1998}, together with the transformations
in~\cite{GI}, reduces the cyclicity problem for the period
annulus to determining the maximal number of zeros of the
following generating function:
\begin{equation}\label{eq:I}
\begin{split}
I(h)=\mu_1hI_{0,0}(h)+\mu_2I_{1,0}(h)
           +\mu_3I_{0,1}(h)
+\mu_4\{2I_{-1,0}(h)+3\kappa hI_{-1,1}(h)\},
\end{split}
\end{equation}
where $\mu_1,\ldots,\mu_4$ are arbitrary real constants,
as in~\cite[Eq.~(9)]{GI} and~\cite[Eq.~(4)]{Zhao}.
Zeros are counted with multiplicity, and the identically
zero function is excluded.

This reduction applies to analytic quadratic perturbations
even when their first-order Melnikov function vanishes
identically. The coefficients $\mu_1,\ldots,\mu_4$
in~\eqref{eq:I} can be prescribed independently through
suitable quadratic perturbations. Simple zeros of $I(h)$ then
yield hyperbolic limit cycles.

The upper bound on cyclicity will follow from the following
theorem. Its proof is completed in
Section~\ref{subsec:generating-proof}.

\begin{theorem}\label{thm:generating}
For every fixed $\kappa>1$, if the function $I(h)$
in~\eqref{eq:I} is not identically zero, then it has at most
three zeros on $(h_c,h_s)$, counted with multiplicity.
\end{theorem}

\subsection{Reduction to linear combinations of four functions}

To prove Theorem~\ref{thm:generating}, we shall estimate the
number of zeros of $I(h)$ by studying the zeros of an auxiliary
function $\overline G(h)$. Following~\cite{GI,Zhao}, we use
the identity
\begin{equation}\label{eq:RolleI}
\frac d{dh}\left(\frac{I(h)}h\right)
       =\frac{\overline G(h)}{h^2},
\end{equation}
where
\begin{equation}\label{eq:G}
\begin{split}
\overline G(h)={}&(\widetilde\mu_1h^2+\widetilde\mu_3)I'_{0,0}(h)
                   +\widetilde\mu_2I'_{1,1}(h)\\
&+\widetilde\mu_4\{-4hI'_{-1,0}(h)
                   +(3\kappa h^2-4)I'_{-1,1}(h)\}.
\end{split}
\end{equation}
Here primes denote derivatives with respect to $h$,
and the coefficients are related by the invertible change
\begin{equation*}
\begin{split}
\widetilde\mu_1=\mu_1,\,\,\,\,
\widetilde\mu_2=-\frac23\mu_2-
             \frac{2(\kappa-1)}{3\kappa}\mu_3,
\,\,\,\, \widetilde\mu_3=-\frac{2}{3\kappa}\mu_3,\,\,\,\, \widetilde\mu_4=\mu_4.
\end{split}
\end{equation*}
With this coefficient redefinition, the expression for
$\overline G$ agrees with that given after Eq.~(9)
in~\cite{GI}. The change is written explicitly after
Eq.~(20) in~\cite{Zhao}.
We retain the tildes to distinguish the coefficients of
$\overline G(h)$ from those of $I(h)$.

Following~\cite{GI,Zhao}, we introduce a new variable $s$ by
\begin{equation}\label{eq:s-variable}
s=\frac{9\kappa}{4}h^2,\qquad
h=h(s)=-\frac23\sqrt{\frac{s}{\kappa}}.
\end{equation}
Then $1<s<\kappa$, with $s=\kappa$ corresponding to the
center and $s=1$ to the outer boundary.

We use the notation $J_1,J_2$ of~\cite[p.~74]{GI}
and $J_3,J_4$ of~\cite[Eq.~(18)]{Zhao}:
\begin{equation}\label{eq:periods}
\begin{aligned}
J_1(s)&=I'_{0,0}(h(s)),&
J_2(s)&=I'_{1,1}(h(s)),\\
J_3(s)&=I'_{-1,0}(h(s)),&
J_4(s)&=I'_{-1,1}(h(s)).
\end{aligned}
\end{equation}
The primes denote derivatives with respect to $h$.
Following~\cite[Eq.~(12)]{Zhao}, write
$G(s)=\overline G(h(s))$. Substituting~\eqref{eq:s-variable}
and~\eqref{eq:periods} into~\eqref{eq:G}, we obtain
\begin{equation}\label{eq:G-s}
G(s)={}\left(\frac{4\widetilde\mu_1}{9\kappa}s
             +\widetilde\mu_3\right)J_1(s)
        +\widetilde\mu_2J_2(s)
+\frac43\widetilde\mu_4
  \left\{2\sqrt{\frac{s}{\kappa}}J_3(s)+(s-3)J_4(s)\right\}.
\end{equation}
The integral $J_1(s)$ satisfies~\cite[proof of Proposition~5]{GI}
\begin{equation}\label{eq:J1-positive}
J_1(s)>0,\qquad 1<s<\kappa.
\end{equation}
Keeping the ratio $w$ of~\cite[Eq.~(29)]{Zhao}, we introduce
the auxiliary ratio $v$:
\begin{equation}\label{eq:ratios}
w(s)=\frac{J_2(s)}{J_1(s)},\qquad
v(s)=\frac{2\sqrt{s/\kappa}\,J_3(s)+(s-3)J_4(s)}{J_1(s)}.
\end{equation}
Then~\eqref{eq:G-s} becomes
\begin{equation}\label{eq:G-ratios}
\frac{G(s)}{J_1(s)}
 =\frac{4\widetilde\mu_1}{9\kappa}s+\widetilde\mu_3
  +\widetilde\mu_2w(s)+\frac43\widetilde\mu_4v(s).
\end{equation}

Under the change of variable $h=h(s)$, the zeros of
$\overline G(h)$ correspond to those of $G(s)$ with the same
multiplicities. By~\eqref{eq:J1-positive}, division by $J_1(s)$
preserves the zeros of $G(s)$ and their multiplicities.
Thus, by~\eqref{eq:G-ratios}, the problem reduces to estimating
the number of zeros of linear combinations of
$1$, $s$, $w(s)$ and $v(s)$.

\subsection{The Chebyshev property and the Wronskian inequality}\label{subsec:chebyshev-reduction}

To estimate the number of zeros of the linear combinations
in~\eqref{eq:G-ratios}, we shall prove that the ordered system
$(1,s,w(s),v(s))$ is an extended complete Chebyshev system on
$(1,\kappa)$.

For functions $f_0,\ldots,f_{k-1}$ of $s$, we define their
Wronskian determinant by
\[
W(f_0,\ldots,f_{k-1})(s)
   =\det\bigl(f_j^{(i)}(s)\bigr)_{0\leq i,j\leq k-1}.
\]
Here $f_j^{(i)}$ denotes the $i$-th derivative of $f_j$
with respect to $s$, with $f_j^{(0)}=f_j$.
The four successive Wronskians of
$(1,s,w,v)$ are $W(1)=W(1,s)=1,\,\,W(1,s,w)=d^2w/ds^2$ and
\begin{equation}\label{eq:target}
W(1,s,w,v)=\frac{d^2w}{ds^2}\frac{d^3v}{ds^3}
            -\frac{d^3w}{ds^3}\frac{d^2v}{ds^2}.
\end{equation}
The definition of an ECT-system and the Wronskian criterion
are recalled in Section~\ref{sec:prelim}. By Lemma~\ref{lem:ECT},
it suffices to show that the last two Wronskians do not vanish.
The inequality $d^2w/ds^2<0$ follows from the concavity
recalled in Section~\ref{sec:comparison}. The remaining step is the
following strict Wronskian inequality. Its proof is completed
in Section~\ref{subsec:chebyshev-proof} using the positivity
established in Section~\ref{sec:positive}.

\begin{theorem}\label{thm:Wronskian}
For every $\kappa>1$ and $1<s<\kappa$, we have $W(1,s,w,v)>0$.
\end{theorem}

By~\eqref{eq:G-ratios}, the Chebyshev property implies that
$G(s)$ has at most three zeros, counted with multiplicity,
provided that $G(s)\not\equiv0$. Finally,
Lemma~\ref{lem:I-center} gives $I(h_c)=0$. Since $h\ne0$,
Rolle's theorem applied to~\eqref{eq:RolleI}
transfers this bound from $\overline G(h)$ to $I(h)$. If
$\overline G(h)\equiv0$, the same identity and the endpoint value
force $I(h)\equiv0$. These steps are carried out in
Section~\ref{sec:conclusion}.

The main remaining task is therefore the proof of
Theorem~\ref{thm:Wronskian}. After collecting the necessary
equations in Section~\ref{sec:prelim}, we derive a second-order
inhomogeneous equation and construct a comparison function in
Section~\ref{sec:comparison}. In Section~\ref{sec:positive},
we assume that the Wronskian has a zero and use the
comparison equation and two polynomial sign identities to
obtain a contradiction. The argument shows that the Wronskian
has no zeros on the interval.

\section{Preliminary results}\label{sec:prelim}

We collect the equations, center expansions and results on
zeros used in the proof outlined above.

\subsection{Picard--Fuchs system and the Riccati equation}

For the comparison argument, we introduce the normalized
variable $q$ and parameter $p$ by
\begin{equation}\label{eq:coordinates}
q=\frac{\kappa-s}{\kappa-1},\qquad
p=\frac{\kappa-1}{\kappa},\qquad 1<s<\kappa,\quad \kappa>1.
\end{equation}
The inverse change is
\begin{equation}\label{eq:inverse-coordinates}
\kappa=\frac1{1-p},\qquad
s=s(q)=\frac{1-pq}{1-p},\qquad 0<p<1,\quad 0<q<1.
\end{equation}
The center corresponds to $q=0$,
and the outer boundary to $q=1$.
The center expansions in~\cite[Lemma~8]{Zhao}, derived from
the Picard--Fuchs identities in~\cite[p.~72]{GI}, give the
common positive center limit
\[
J_i(\kappa)=J_1(\kappa)>0,
\qquad i=1,\ldots,4.
\]
Define the normalized integrals
\begin{equation}\label{eq:normalized-integrals}
\mathcal J_i(q)=\frac{J_i(s(q))}{J_1(\kappa)},\qquad
 i=1,2,3,4.
\end{equation}
They satisfy
\begin{equation}\label{eq:normalized-center-values}
\mathcal J_i(0)=1,\qquad i=1,2,3,4.
\end{equation}
In the normalized proof below, $p\in(0,1)$ is fixed and
$q\in(0,1)$ is the independent variable. We return to the
original variables in Section~\ref{sec:conclusion}.
The ratios in~\eqref{eq:ratios} become
\begin{equation}\label{eq:normalized-ratios}
w(q)=\frac{\mathcal J_2(q)}{\mathcal J_1(q)},\qquad
v(q)=\frac{2(1-p)\sqrt{1-pq}\,\mathcal J_3(q)
                 +(-2+3p-pq)\mathcal J_4(q)}{(1-p)\mathcal J_1(q)}.
\end{equation}
The first four equations of~\cite[Eq.~(31)]{Zhao} form
a closed system, which can be written as
\begin{equation}\label{eq:PF-original}
(9h^2-4)(9\kappa h^2-4)
\begin{pmatrix}
I''_{0,0}\\ I''_{1,1}\\ I''_{-1,0}\\ I''_{-1,1}
\end{pmatrix}
=M_4(h)
\begin{pmatrix}
I'_{0,0}\\ I'_{1,1}\\ I'_{-1,0}\\ I'_{-1,1}
\end{pmatrix},
\end{equation}
where
\[
M_4(h)=
\begin{pmatrix}
-3h(9\kappa h^2-4)&12(\kappa-1)h&0&0\\
-3h(9\kappa h^2-4)&3h(9\kappa h^2-4)&0&0\\
8(\kappa-1)&-8(\kappa-1)&-6\kappa h(9h^2-4)&2\kappa(9h^2-4)\\
-12(\kappa-1)h&12(\kappa-1)h&4(9h^2-4)&-3\kappa h(9h^2-4)
\end{pmatrix}.
\]
Primes denote derivatives with respect to $h$.

The changes~\eqref{eq:s-variable} and~\eqref{eq:coordinates}
express $h$ as the composite function $h(s(q))$.
Together with the definitions~\eqref{eq:periods} and
\eqref{eq:normalized-integrals}, these changes transform
\eqref{eq:PF-original} into the Picard--Fuchs system
\begin{equation}\label{eq:PF-normalized}
6q(1-q)
\begin{pmatrix}
\mathcal J_1'\\ \mathcal J_2'\\ \mathcal J_3'\\ \mathcal J_4'
\end{pmatrix}
=
\begin{pmatrix}
q-1&1&0&0\\
q-1&1-q&0&0\\
-\dfrac1{\sqrt{1-pq}}&\dfrac1{\sqrt{1-pq}}&2q&\dfrac q{\sqrt{1-pq}}\\[2pt]
-1&1&\dfrac{2(1-p)q}{\sqrt{1-pq}}&q
\end{pmatrix}
\begin{pmatrix}
\mathcal J_1\\ \mathcal J_2\\ \mathcal J_3\\ \mathcal J_4
\end{pmatrix}.
\end{equation}
Primes denote derivatives with respect to $q$.
In particular, $1-pq>0$ on the parameter domain.

Differentiating the formula for $w$ in~\eqref{eq:normalized-ratios}
and using the first two equations of~\eqref{eq:PF-normalized}, we obtain
\begin{equation}\label{eq:Riccati}
6q(1-q)w'=q-1+2(1-q)w-w^2.
\end{equation}
When~\eqref{eq:Riccati} is solved for $w'$, its apparent
singularity at $q=0$ is removable along the analytic center solution. The Riccati
equation will be used to eliminate derivatives of $w$ in the
sign computations of Section~\ref{sec:positive}.

\subsection{Series expansions at the center}

Since $q=0$ corresponds to the nondegenerate center, the four
normalized integrals $\mathcal J_i(q)$ are analytic there, with
the center values given in~\eqref{eq:normalized-center-values}.
Eliminating $\mathcal J_2$ from the first two equations
of~\eqref{eq:PF-normalized} gives
\begin{equation}\label{eq:Jode}
q(1-q)\mathcal J_1''+(1-2q)\mathcal J_1'-\frac5{36}\mathcal J_1=0.
\end{equation}
Substitution of a power series into~\eqref{eq:Jode} gives
\begin{equation}\label{eq:Jseries}
\mathcal J_1(q)=\sum_{n=0}^{\infty}j_nq^n,\,\,\,\,
j_0=1,\,\,\,
j_{n+1}=\frac{(6n+1)(6n+5)}{36(n+1)^2}j_n,\,\,\, n=0,1,\ldots.
\end{equation}
The recurrence gives $j_n>0$ for every $n\geq0$ and
$j_{n+1}/j_n\to1$, which show that  the series converges for $|q|<1$.
Its sum agrees with the center solution near zero and, by
uniqueness for~\eqref{eq:Jode} on $0<q<1$, represents
$\mathcal J_1(q)$ throughout $0\leq q<1$.
The first equation of~\eqref{eq:PF-normalized} and the
recurrence in~\eqref{eq:Jseries} yield
\begin{equation}\label{eq:J2series}
\mathcal J_2(q)=(1-q)\bigl(\mathcal J_1(q)+6q\mathcal J_1'(q)\bigr)=1-\sum_{n=1}^{\infty}\frac{j_n}{6n-1}q^n.
\end{equation}
Using~\eqref{eq:Jseries}, \eqref{eq:J2series} and the formula
for $w$ in~\eqref{eq:normalized-ratios}, we obtain the center expansions
\begin{equation}\label{eq:centerJE}
\begin{array}{cl}
\mathcal J_1(q)&=1+\dfrac5{36}q+\dfrac{385}{5184}q^2+O(q^3),
\\[2ex]
\mathcal J_2(q)&=1-\dfrac1{36}q-\dfrac{35}{5184}q^2+O(q^3),
\end{array}
\end{equation}
and
\begin{equation}\label{eq:centerw}
w(q)=1-\frac q6-\frac{25}{432}q^2 -\frac{775}{23328}q^3+O(q^4).
\end{equation}
In particular,
\begin{equation}\label{eq:centerw-derivatives}
w''(0)=-\frac{25}{216},\qquad
w'''(0)=-\frac{775}{3888}.
\end{equation}
To obtain the center expansion of $v$, first note
from~\eqref{eq:centerJE} that
\[
\frac{\mathcal J_2(q)-\mathcal J_1(q)}q=-\frac16+O(q).
\]
For $q>0$, dividing the last two equations
of~\eqref{eq:PF-normalized} by $q$ gives
\[
\begin{aligned}
6(1-q)\mathcal J_3'
 &=\frac{(\mathcal J_2-\mathcal J_1)/q+\mathcal J_4}
         {\sqrt{1-pq}}+2\mathcal J_3,\\
6(1-q)\mathcal J_4'
 &=\frac{\mathcal J_2-\mathcal J_1}{q}
   +\frac{2(1-p)}{\sqrt{1-pq}}\mathcal J_3+\mathcal J_4.
\end{aligned}
\]
Letting $q\to0$ and using~\eqref{eq:normalized-center-values},
we obtain
\[
\begin{aligned}
6\mathcal J_3'(0)=-\frac16+1+2=\frac{17}{6},\,\,\,\,
6\mathcal J_4'(0)=-\frac16+2(1-p)+1=\frac{17-12p}{6}.
\end{aligned}
\]
Thus
\[
\mathcal J_3(q)=1+\frac{17}{36}q+O(q^2),\qquad
\mathcal J_4(q)=1+\frac{17-12p}{36}q+O(q^2).
\]
Using~\eqref{eq:Jseries} and~\eqref{eq:J2series}, we continue
the coefficient comparison in the last two equations
of~\eqref{eq:PF-normalized} to determine the Taylor coefficients
of $\mathcal J_3$ and $\mathcal J_4$ through order $q^3$.
Substituting the Taylor expansions of $\mathcal J_1$,
$\mathcal J_3$ and $\mathcal J_4$ into the formula for $v$
in~\eqref{eq:normalized-ratios}, and expanding $\sqrt{1-pq}$
at $q=0$, we obtain
\begin{equation}\label{eq:centerv}
v(q)=\frac{p}{1-p}(1-q)-\frac{5p}{36}q^2
     -\frac{5p(9p+22)}{1944}q^3+O(q^4).
\end{equation}
The expansions~\eqref{eq:centerw} and~\eqref{eq:centerv}
can also be deduced from~\cite[Lemma~8]{Zhao}. That lemma
gives the Taylor expansions of $J_i(s)$, $i=1,\ldots,4$,
at $s=\kappa$ through degree four in $s-\kappa$.
To obtain the expansions of $w$ and $v$, first replace
$s-\kappa$ by $-(\kappa-1)q$ in the expansions of $J_i(s)$
and divide by $J_1(\kappa)$, according
to~\eqref{eq:coordinates} and~\eqref{eq:normalized-integrals}.
Then substitute the Taylor expansions of
$\mathcal J_1(q),\ldots,\mathcal J_4(q)$ into the formulas
for $w$ and $v$ in~\eqref{eq:normalized-ratios} to
obtain~\eqref{eq:centerw} and~\eqref{eq:centerv}.
For completeness and the reader's convenience, we include a direct
derivation of these Taylor expansions in the present variables.

The explicit recurrence for all coefficients $j_n$
in~\eqref{eq:Jseries} was not given in~\cite{Zhao}.
Here we derive this recurrence from~\eqref{eq:Jode}
and prove that the series $\sum_{n=0}^{\infty}j_nq^n$
converges for $|q|<1$.

\subsection{Chebyshev systems and Rolle's theorem}

\begin{definition}\label{def:Chebyshev}
Let $f_0,\ldots,f_{n-1}$ be real analytic functions on an
interval $(a,b)$. Their ordered tuple is a
\emph{Chebyshev system} if each linear combination with nonzero
coefficient vector has at most $n-1$ distinct zeros on
$(a,b)$. It is a \emph{complete Chebyshev system} if, for each
$k=1,\ldots,n$, every linear combination of
$f_0,\ldots,f_{k-1}$ with nonzero coefficient vector has at most
$k-1$ distinct zeros on $(a,b)$. It is an
\emph{extended complete Chebyshev system} (ECT-system) if these
bounds remain valid when zeros are counted with multiplicity.
\end{definition}

We use the Wronskian notation introduced in
Section~\ref{subsec:chebyshev-reduction}, now with derivatives
taken with respect to $q$. The following standard criterion
is given in~\cite[Lemma~2.3(b)]{GMV}.

\begin{lemma}[Wronskian criterion]\label{lem:ECT}
The ordered analytic tuple $(f_0,\ldots,f_{n-1})$ is an
ECT-system on $(a,b)$ if and only if
\[
W(f_0,\ldots,f_{k-1})(q)\ne0,
\qquad q\in(a,b),\quad k=1,\ldots,n.
\]
\end{lemma}

We shall also use the following direct consequence of Rolle's
theorem.

\begin{lemma}\label{lem:endpoint}
Let $F$ be continuous on $[a,b)$ and real analytic on $(a,b)$,
with $F(a)=0$ and $F'\not\equiv0$. If $F'$ has at most $m$
zeros on $(a,b)$, counted with multiplicity, then $F$ has at
most $m$ zeros there, counted with multiplicity.
\end{lemma}

We recall the following standard result on functions with
prescribed zeros; see~\cite[Definition~2.1
and Proposition~2.2, p.~178]{Mazure2005}.
For completeness, we give below a direct algebraic proof
that does not use the results of that reference.

\begin{proposition}\label{prop:prescribed-zeros}
Let $V$ be an $n$-dimensional real vector space of real analytic
functions on $(a,b)$, where $n\geq2$. Suppose that every nonzero
function in $V$ has at most $n-1$ zeros on $(a,b)$, counted
with multiplicity. Then, for any prescribed points
$a<t_1<\cdots<t_{n-1}<b$, there is a nonzero function in $V$
whose zeros on $(a,b)$ are precisely these points, all simple.
\end{proposition}

\begin{proof}
Let
\[
L:V\longrightarrow\mathbb R^{n-1}
\]
be the linear map defined by
\[
L(f)=\bigl(f(t_1),\ldots,f(t_{n-1})\bigr),
\qquad f\in V.
\]
Clearly, $L$ is linear. Since $\dim V=n$ and
$\dim\mathbb R^{n-1}=n-1$, we have
\(\operatorname{rank}L\leq n-1.\)
By the rank-nullity theorem,
\(\dim V=\dim\ker L+\operatorname{rank}L.\)
Consequently,
\[
\dim\ker L
=n-\operatorname{rank}L\geq1.
\]
Hence there exists a nonzero function $f\in V$ such that
\[
f(t_1)=\cdots=f(t_{n-1})=0.
\]

We now show that these are precisely the zeros of $f$
in $(a,b)$ and that they are all simple. Indeed, since
$f\not\equiv0$ is real analytic, each of its zeros has
a well-defined finite multiplicity. By construction,
$f$ has at least $n-1$ zeros in $(a,b)$, counted with
multiplicity. If any one of the prescribed zeros had
multiplicity greater than one, or if $f$ had an additional
zero in $(a,b)$, then $f$ would have at least $n$ zeros
in $(a,b)$, counted with multiplicity. This contradicts
the hypothesis that every nonzero element of $V$ has
at most $n-1$ zeros there.
Therefore, the zeros of $f$ in $(a,b)$ are precisely
$t_1,\ldots,t_{n-1}$, and all of them are simple.
\end{proof}

\section{A forced equation and its comparison function}\label{sec:comparison}

In this section we introduce an auxiliary function $U$
to study the Wronskian in~\eqref{eq:target}.
We derive a second-order inhomogeneous equation for $U$ and
construct a comparison function. We then recall the concavity
of $w$ and establish the further sign properties needed for
the Wronskian argument in Section~\ref{sec:positive}.

\subsection{A second-order equation and a comparison function}

For the function $v(q)$ in~\eqref{eq:normalized-ratios},
we subtract its first-degree Taylor polynomial at the center
$q=0$ and multiply the difference by $-1/p$. This defines
\begin{equation}\label{eq:U}
U(q)=\frac{1-q}{1-p}-\frac{v(q)}p,
\end{equation}
so that
\begin{equation}\label{eq:vU}
v=\frac{p}{1-p}(1-q)-pU.
\end{equation}
By~\eqref{eq:centerv}, we have $U(q)=O(q^2)$, and hence
\begin{equation}\label{eq:U-initial}
U(0)=U'(0)=0.
\end{equation}
The affine term in~\eqref{eq:vU} disappears after two
differentiations.

By~\eqref{eq:inverse-coordinates},
\[
\frac d{ds}=-\frac{1-p}{p}\frac d{dq}.
\]
Differentiating~\eqref{eq:vU} twice and three times with
respect to $q$ gives $v''=-pU''$ and $v'''=-pU'''$.
Consequently,
\begin{equation}\label{eq:conversion}
W(1,s,w)=\frac{d^2w}{ds^2}
 =\frac{(1-p)^2}{p^2}w'',\,\,\,
W(1,s,w,v)
 =\frac{(1-p)^5}{p^4}
   \bigl(w''U'''-w'''U''\bigr).
\end{equation}
The Wronskians in these identities are computed with respect to $s$.
Since $0<p<1$ and $W(1)=W(1,s)=1$, it follows
from~\eqref{eq:conversion} and Lemma~\ref{lem:ECT} that,
to prove the ECT property of $(1,s,w,v)$, it suffices to prove
$w''<0$ and $w''U'''-w'''U''>0$ on $(0,1)$.

The operator in the following proposition is obtained from
the $s$-form of the second-order operator of Gavrilov and
Iliev~\cite[p.~74]{GI} by~\eqref{eq:coordinates} and multiplication
by $-p$; see also~\cite[Eq.~(13)]{Zhao}.
We use this operator to derive equations for $\mathcal J_1U$
and $\mathcal J_1U_*$, where $U_*$ is a comparison function
constructed from $w$ in the following proposition.
Expanding the product derivatives and dividing by
$\mathcal J_1>0$ gives equations for $U$ and $U_*$
with the same differential operator, as displayed
in~\eqref{eq:pair} below. These comparison equations will be used
to prove $w''U'''-w'''U''>0$ in Section~\ref{sec:positive}.

\begin{proposition}\label{prop:forced}
Let
\begin{equation}\label{eq:Bp}
\mathcal B_p=(1-pq)(1-q)\frac{d^2}{dq^2}
             -\frac{1-p}{2}\frac d{dq}+\frac{5p}{36}.
\end{equation}
Then
\begin{equation}\label{eq:Zeq}
\mathcal B_p(\mathcal J_1U)=\frac5{18}\mathcal J_1.
\end{equation}
Moreover, define the comparison function
\begin{equation}\label{eq:Ustar}
U_*(q)=\frac{12}{5}\left(1-w(q)-\frac q6\right).
\end{equation}
It satisfies
\begin{equation}\label{eq:Zstareq}
\mathcal B_p(\mathcal J_1U_*)=\frac5{18}\mathcal J_1+(1-p)q\mathcal J_1',
\end{equation}
with initial conditions
\begin{equation}\label{eq:Ustar-initial}
U_*(0)=U_*'(0)=0.
\end{equation}
\end{proposition}
\begin{proof}
We recall the equation proved in~\cite[p.~72]{GI} for
\[
K(h)=-4h\frac{dI_{-1,0}}{dh}
       +(3\kappa h^2-4)\frac{dI_{-1,1}}{dh}.
\]
It is
\begin{equation}\label{eq:GI-K-equation}
L_2(h)K=\frac43(\kappa-1)
\left[h(9\kappa h^2-4)\frac{d^3I_{1,1}}{dh^3}
      +(6\kappa h^2+8)\frac{d^2I_{1,1}}{dh^2}\right],
\end{equation}
where
\[
L_2(h)=h(9\kappa h^2-4)\frac{d^2}{dh^2}
       -(9\kappa h^2-8)\frac d{dh}+5\kappa h.
\]
Under~\eqref{eq:s-variable}, the definitions
\eqref{eq:periods} and~\eqref{eq:ratios} give
\[
K=\frac43J_1v,\,\,\,\,L_2(h)=-36\kappa hL_2,\] where
\[
L_2=s(1-s)\frac{d^2}{ds^2}-\frac12\frac d{ds}-\frac5{36}.
\]
Here and below, $L_2$ acts on functions of $s$, and
$\mathcal B_p$ acts on functions of $q$.
Consequently, equation~\eqref{eq:GI-K-equation} becomes
\begin{equation}\label{eq:GI-s-equation}
L_2(J_1v)=(\kappa-1)
\left[s(1-s)\frac{d^2J_2}{ds^2}
      -\frac{5s+3}{6}\frac{dJ_2}{ds}\right].
\end{equation}
The first two equations of~\eqref{eq:PF-original}, written
in the variable $s$ using~\eqref{eq:s-variable}, give
\begin{equation}\label{eq:PF-s-first-relation}
(\kappa-1)\frac{dJ_2}{ds}
=(s-1)\frac{dJ_1}{ds}+\frac16J_1.
\end{equation}
Differentiating this relation with respect to $s$ gives
\begin{equation}\label{eq:PF-s-second-relation}
(\kappa-1)\frac{d^2J_2}{ds^2}
=(s-1)\frac{d^2J_1}{ds^2}+\frac76\frac{dJ_1}{ds}.
\end{equation}
Substituting~\eqref{eq:PF-s-first-relation}
and~\eqref{eq:PF-s-second-relation}
into~\eqref{eq:GI-s-equation} yields
\begin{equation}\label{eq:L2-difference}
L_2(J_1v)-L_2\bigl((s-1)J_1\bigr)=\frac5{18}J_1.
\end{equation}
Finally, for any twice differentiable function $f(s)$,
the affine change~\eqref{eq:inverse-coordinates} with $p$ fixed
and the chain rule give
\[
\frac d{dq}f(s(q))=-\frac p{1-p}\frac{df}{ds}(s(q)),\,\,\,
\frac{d^2}{dq^2}f(s(q))
 =\frac{p^2}{(1-p)^2}\frac{d^2f}{ds^2}(s(q)).
\]
Together with
\[
s(q)(1-s(q))=-\frac{p(1-pq)(1-q)}{(1-p)^2},
\]
these identities and~\eqref{eq:Bp} yield
\begin{equation}\label{eq:operator-change}
\mathcal B_p\bigl[f(s(q))\bigr]
=-p\,(L_2f)(s(q)).
\end{equation}
To express the first term in~\eqref{eq:U} in terms of $s(q)$,
we use~\eqref{eq:inverse-coordinates} to write
\[
\frac{1-q}{1-p}=\frac{s(q)-1}{p}.
\]
Substituting the above equality into~\eqref{eq:U} gives
\[
U(q)=\frac{s(q)-1-v(s(q))}{p}.
\]
Multiplying by $\mathcal J_1(q)$ and
using~\eqref{eq:normalized-integrals}, we obtain
\[
\mathcal J_1(q)U(q)
=\frac{\bigl[(s-1)J_1(s)-J_1(s)v(s)\bigr]_{s=s(q)}}
       {pJ_1(\kappa)}.
\]
Applying~\eqref{eq:operator-change} and then
using~\eqref{eq:L2-difference}, we obtain
\[
\mathcal B_p(\mathcal J_1U)(q)
 =\frac{\bigl[L_2(J_1v)-L_2((s-1)J_1)\bigr]_{s=s(q)}}
          {J_1(\kappa)}
 =\frac5{18}\frac{J_1(s(q))}{J_1(\kappa)}
 =\frac5{18}\mathcal J_1(q).
\]
The operators inside the bracket act with respect to $s$
before evaluation at $s=s(q)$.
This proves~\eqref{eq:Zeq}.

To prove~\eqref{eq:Zstareq}, we first use~\eqref{eq:Ustar}
and~\eqref{eq:normalized-ratios} to obtain
\begin{equation}\label{eq:J1Ustar-product}
\mathcal J_1U_* =\frac{12}{5}
 \left(\mathcal J_1-\mathcal J_2-\frac q6\mathcal J_1\right).
\end{equation}
Subtracting the second equation of~\eqref{eq:PF-normalized}
from the first gives the first relation below, while rewriting
the first equation gives the second:
\[
\mathcal J_1'-\mathcal J_2'
=\frac{\mathcal J_2}{6(1-q)},\qquad
\frac{\mathcal J_2}{1-q}
=\mathcal J_1+6q\mathcal J_1'.
\]
Differentiating~\eqref{eq:J1Ustar-product} and using these two
identities, we obtain
\[(\mathcal J_1U_*)'=2q\mathcal J_1'.\]
Combining this identity with~\eqref{eq:Jode} yields
\[
(\mathcal J_1U_*)''=\frac{2q}{1-q}\mathcal J_1'
                    +\frac{5}{18(1-q)}\mathcal J_1.
\]
Substituting these derivative formulas and~\eqref{eq:J2series}
into~\eqref{eq:Bp} proves~\eqref{eq:Zstareq}.
The initial conditions~\eqref{eq:Ustar-initial} follow from
the center expansion~\eqref{eq:centerw}.
\end{proof}

\subsection{Concavity and monotonicity of the first ratio}

To express the coefficients of the comparison equations, define
the logarithmic derivative
\begin{equation}\label{eq:eta}
\eta(q)=\frac{\mathcal J_1'(q)}{\mathcal J_1(q)}.
\end{equation}
It is well defined for $0\leq q<1$ since $\mathcal J_1>0$.
The first equation of~\eqref{eq:PF-normalized} yields
\begin{equation}\label{eq:eta-PF}
\eta(q)=\frac{q+w-1}{6q(1-q)}.
\end{equation}

\begin{lemma}\label{lem:ratio}
On $0<q<1$,
\begin{equation}\label{eq:wrange}
1-q<w(q)<1,\qquad w''(q)<0,\qquad
\frac d{dq}(q\eta(q))>0.
\end{equation}
\end{lemma}
\begin{proof}
For $0<q<1$, the positive coefficients in~\eqref{eq:Jseries}
give $\mathcal J_1(q)>1$ and $\mathcal J_1'(q)>0$.
Hence $\eta(q)>0$ by~\eqref{eq:eta}, which implies
$w(q)>1-q$ by~\eqref{eq:eta-PF}.
Equation~\eqref{eq:J2series} gives $\mathcal J_2(q)<1$, so
$w(q)=\mathcal J_2(q)/\mathcal J_1(q)<1$.
Together these observations prove the first inequality
in~\eqref{eq:wrange}.

For concavity, we use~\cite[Lemma~15]{Zhao}, where the ratio
$w(s)=J_2(s)/J_1(s)$ has independent variable $s$.
Here we use the same symbol for the composition
$w(q)=w(s(q))$. Since $s(q)$ is affine, the chain rule
and~\eqref{eq:inverse-coordinates} give
\[
\frac d{dq}w(s(q))
 =-\frac p{1-p}\left.\frac{dw}{ds}\right|_{s=s(q)},\,\,\,
\frac{d^2}{dq^2}w(s(q))
 =\frac{p^2}{(1-p)^2}
   \left.\frac{d^2w}{ds^2}\right|_{s=s(q)}.
\]
Thus the first derivative changes sign, while the second
derivative retains its sign. The inequality $d^2w/ds^2<0$
proved in~\cite[Lemma~15]{Zhao} therefore gives $w''(q)<0$,
which is the second inequality in~\eqref{eq:wrange}.

Using~\eqref{eq:Jseries} and $j_0j_1>0$, we obtain by
differentiation and symmetrization
\[
\frac d{dq}(q\eta(q))
 =\frac{1}{2q\mathcal J_1(q)^2}
   \sum_{m,n\geq0}(m-n)^2j_mj_nq^{m+n}>0,\qquad 0<q<1.\qedhere
\]
\end{proof}

Put
\begin{equation}\label{eq:abc}
a=(1-pq)(1-q),\,\,
b=2a\eta-\frac{1-p}{2},\,\,
c=a(\eta'+\eta^2)-\frac{1-p}{2}\eta+\frac{5p}{36}.
\end{equation}
By~\eqref{eq:eta}, we have
$\mathcal J_1''/\mathcal J_1=\eta'+\eta^2$.
The product rule therefore gives
\[
\frac{(\mathcal J_1U)'}{\mathcal J_1}=U'+\eta U,\,\,\,
\frac{(\mathcal J_1U)''}{\mathcal J_1}
 =U''+2\eta U'+(\eta'+\eta^2)U.
\]
Substituting into~\eqref{eq:Bp} yields
$\mathcal B_p(\mathcal J_1U)/\mathcal J_1=aU''+bU'+cU$,
with the coefficients in~\eqref{eq:abc}.
The same calculation applies to $U_*$.
Dividing~\eqref{eq:Zeq} and~\eqref{eq:Zstareq}
by $\mathcal J_1$ now gives the scalar comparison equations
\begin{equation}\label{eq:pair}
aU''+bU'+cU=\frac5{18},\qquad
aU_*''+bU_*'+cU_*=\frac5{18}+(1-p)q\eta.
\end{equation}
All coefficients are analytic on $0\leq q<1$, and $a>0$.
By~\eqref{eq:Ustar}, Lemma~\ref{lem:ratio} and the initial
conditions~\eqref{eq:Ustar-initial}, we have, for $0<q<1$,
\begin{equation}\label{eq:positivefacts}
\begin{aligned}
U_*''(q)&=-\frac{12}{5}w''(q)>0,\\
U_*'(q)&=\int_0^q U_*''(t)\,dt>0,\\
U_*(q)&=\int_0^q(q-t)U_*''(t)\,dt>0.
\end{aligned}
\end{equation}

\section{Positivity of the Wronskian}\label{sec:positive}

In this section we prove that $w''U'''-w'''U''>0$ on $(0,1)$.
We first establish two algebraic inequalities using the
Riccati equation. We then combine them with the center
expansions and the comparison function $U_*$ to rule out
a first zero of the Wronskian. By~\eqref{eq:conversion},
this gives $W(1,s,w,v)>0$.

\subsection{Two algebraic sign identities}

Throughout this section primes are $q$-derivatives with $p$
fixed, and $\eta,w$ satisfy~\eqref{eq:eta-PF} and~\eqref{eq:Riccati}.
\begin{lemma}\label{lem:signs}
For $0<p<1$ and $0<q<1$,
\begin{equation}\label{eq:K-definitions}
\begin{aligned}
K_1&:=U_*\bigl(c(q\eta)'-c'q\eta\bigr)
       -\frac5{18}(q\eta)'<0,\\
K_2&:=aq\eta U_*'''+\bigl((a'+b)q\eta-a(q\eta)'\bigr)U_*''<0.
\end{aligned}
\end{equation}
\end{lemma}
\begin{proof}
Substituting~\eqref{eq:abc}, \eqref{eq:Ustar}
and~\eqref{eq:eta-PF} into the definitions~\eqref{eq:K-definitions}
and eliminating derivatives with~\eqref{eq:Riccati}, we obtain
\begin{equation}\label{eq:Kfactors}
K_1=-\frac{(q+w-1)^2P(p,q,w)}{540q^3(1-q)^2},\,\,\,\,
K_2=-\frac{(q+w-1)^2Q(p,q,w)}{90q^4(1-q)^3}.
\end{equation}
By Lemma~\ref{lem:ratio}, we have $1-q<w<1$.
Since $0<p<1$ and $0<q<1$, all factors in the following
decompositions of the polynomials $P=P(p,q,w)$ and $Q=Q(p,q,w)$
are positive:
\begin{equation*}
\begin{aligned}
P={}&(1-q)\bigl[18q+72(1-w)\bigr]\\
&+q\bigl\{(1-p)[22q+23(1-w)]+31p(q+w-1)\bigr\}>0,
\end{aligned}
\end{equation*}
and
\begin{equation*}
\begin{aligned}
Q={}&(1-p)\bigl\{(q+w-1)[q+2(1-w)]+(1-q)[2q^2+7q(1-w)+(1-w)^2]\\
&\qquad +(1-q)^2(1-w)(4q+1-w)\bigr\}+3p(1-q)\bigl\{(q+w-1)(q+1-w)\\
&\qquad +(1-q)(1-w)(4q+1-w)\bigr\}>0.
\end{aligned}
\end{equation*}
Thus~\eqref{eq:Kfactors} proves both inequalities
in~\eqref{eq:K-definitions}.
\end{proof}

\subsection{The first-zero argument}

Define the $q$-Wronskian
\begin{equation}\label{eq:Wq}
\mathscr W(q)=w''U'''-w'''U''.
\end{equation}

\begin{proposition}\label{prop:global}
For every $0<p<1$,
$\mathscr W(q)>0$ throughout $0<q<1$.
\end{proposition}
\begin{proof}
The expansion~\eqref{eq:centerv}, together with the
definition~\eqref{eq:U}, gives
\[
U''(0)=\frac5{18},\qquad
U'''(0)=\frac{5(9p+22)}{324}.
\]
Combining these values with~\eqref{eq:centerw-derivatives}, we obtain
\begin{equation}\label{eq:Wcenter}
\mathscr W(0)=\frac{125(1-p)}{7776}>0.
\end{equation}
Let
\begin{equation}\label{eq:theta}
\theta(q)=\frac{U''(q)}{U_*''(q)}.
\end{equation}
By~\eqref{eq:positivefacts}, this ratio is well defined.
Moreover,~\eqref{eq:Ustar} and~\eqref{eq:centerw-derivatives} give
\[
U_*''=-\frac{12}{5}w'',\,\,\,
U_*'''=-\frac{12}{5}w''',\,\,\,
U_*''(0)=-\frac{12}{5}w''(0)=\frac5{18}=U''(0).
\]
Using this equality and the two derivative identities above,
together with the derivative of~\eqref{eq:theta}, we obtain
\begin{equation}\label{eq:thetaPrime}
\theta(0)=1,\qquad
\theta'=-\frac{12\mathscr W}{5(U_*'')^2}.
\end{equation}

Suppose, for a contradiction, that $\mathscr W(q)$ has a zero
in $(0,1)$, and let $q_0$ be its first zero. Such a first
zero exists if any zero exists, because of
\eqref{eq:Wcenter} and continuity on every compact subinterval
of $[0,1)$. Since $\mathscr W(0)>0$, continuity and the definition
of $q_0$ give $\mathscr W(q)>0$ for $0<q<q_0$.
By~\eqref{eq:thetaPrime}, we therefore have $\theta'(q)<0$
on $(0,q_0)$ and $\theta'(q_0)=0$.
At $q=q_0$ set
\begin{equation}\label{eq:comparison-ratios}
\alpha=\frac{U(q_0)}{U_*(q_0)},\qquad
\beta=\frac{U'(q_0)}{U_*'(q_0)},\qquad
\delta=\frac{U''(q_0)}{U_*''(q_0)}=\theta(q_0).
\end{equation}
All three denominators are positive by~\eqref{eq:positivefacts}.
We claim
\begin{equation}\label{eq:ordering}
1>\alpha>\beta>\delta.
\end{equation}
Indeed, the initial conditions~\eqref{eq:U-initial}
and~\eqref{eq:Ustar-initial} imply
\[
\alpha=\frac{\displaystyle\int_0^{q_0}(q_0-u)U_*''(u)\theta(u)\,du}
          {\displaystyle\int_0^{q_0}(q_0-u)U_*''(u)\,du},\,\,\,
\beta=\frac{\displaystyle\int_0^{q_0}U_*''(t)\theta(t)\,dt}
          {\displaystyle\int_0^{q_0}U_*''(t)\,dt}.
\]
Since $\theta'(q)<0$ on $(0,q_0)$ and $\theta(0)=1$,
we have $\theta(q_0)<\theta(t)<1$ for $0<t<q_0$.
Using $U_*''>0$ from~\eqref{eq:positivefacts}, we conclude
that the first average $\alpha$ is less than $1$
and the second average $\beta$ exceeds $\delta=\theta(q_0)$.

Finally, write
\begin{equation}\label{eq:integrand-expansion}
(u-t)\bigl[\theta(t)-\theta(u)\bigr]
=(q_0-t)\theta(t)+(q_0-u)\theta(u)
-(q_0-t)\theta(u)-(q_0-u)\theta(t).
\end{equation}
By~\eqref{eq:theta}, we have
\[
\theta(t)U_*''(t)=U''(t),\qquad
\theta(u)U_*''(u)=U''(u).
\]
Multiplying~\eqref{eq:integrand-expansion} by $U_*''(t)U_*''(u)$,
substituting these two identities and integrating over
$0<t<q_0$, $0<u<q_0$, we obtain
\begin{equation}\label{eq:integrated-expansion}
\begin{aligned}
&\int_0^{q_0}\!\int_0^{q_0}
 (u-t)\bigl[\theta(t)-\theta(u)\bigr]
 U_*''(t)U_*''(u)\,dt\,du\\
&=\int_0^{q_0}\!\int_0^{q_0}
 \Bigl[(q_0-t)U''(t)U_*''(u)
       +(q_0-u)U_*''(t)U''(u)\\
&\hspace{18mm}-(q_0-t)U_*''(t)U''(u)
       -(q_0-u)U''(t)U_*''(u)\Bigr]\,dt\,du.
\end{aligned}
\end{equation}
For the term $(q_0-t)U''(t)U_*''(u)$ in
\eqref{eq:integrated-expansion}, we separate the variables
and use the initial conditions~\eqref{eq:U-initial}
and~\eqref{eq:Ustar-initial} to obtain
\[
\begin{aligned}
&\int_0^{q_0}\!\int_0^{q_0}
 (q_0-t)U''(t)U_*''(u)\,dt\,du\\
&\quad=\left(\int_0^{q_0}(q_0-t)U''(t)\,dt\right)
       \left(\int_0^{q_0}U_*''(u)\,du\right)=U(q_0)U_*'(q_0).
\end{aligned}
\]
The double integral of the second term
in~\eqref{eq:integrated-expansion} also equals $U(q_0)U_*'(q_0)$.
The double integral of each of the last two terms equals
$-U_*(q_0)U'(q_0)$.
Consequently,
\begin{equation}\label{eq:alpha-beta-integral}
\begin{aligned}
&\int_0^{q_0}\!\int_0^{q_0}
 (u-t)\bigl[\theta(t)-\theta(u)\bigr]
 U_*''(t)U_*''(u)\,dt\,du\\
&\quad=2\bigl[U_*'(q_0)U(q_0)-U_*(q_0)U'(q_0)\bigr]\\
&\quad=2U_*(q_0)U_*'(q_0)(\alpha-\beta).
\end{aligned}
\end{equation}
Since $\theta$ is strictly decreasing on $(0,q_0)$,
we have $\theta(t)>\theta(u)$ if $0<t<u<q_0$,
and $\theta(t)<\theta(u)$ if $0<u<t<q_0$.
Thus $(u-t)[\theta(t)-\theta(u)]>0$ whenever $t\ne u$.
Since $U_*''(t)U_*''(u)>0$ by~\eqref{eq:positivefacts},
the integrand on the left-hand side of~\eqref{eq:alpha-beta-integral}
is strictly positive off the diagonal $t=u$.
The double integral in~\eqref{eq:alpha-beta-integral}
is therefore strictly positive. By~\eqref{eq:positivefacts},
$U_*(q_0)U_*'(q_0)>0$, and hence $\alpha>\beta$.
This proves \eqref{eq:ordering}, without assuming $\delta\geq0$.

At the first zero $q_0$, \eqref{eq:theta} and
\eqref{eq:comparison-ratios} give
$U'''(q_0)=\delta U_*'''(q_0)$, because $\theta'(q_0)=0$.
Differentiating the first equation in~\eqref{eq:pair} for $U$
and collecting terms gives
\begin{equation}\label{eq:U-equation-derivative}
aU'''+(a'+b)U''+(b'+c)U'+c'U=0.
\end{equation}
At the first zero $q_0$, substituting the relations
in~\eqref{eq:comparison-ratios} and
$U'''(q_0)=\delta U_*'''(q_0)$ into the first equation
in~\eqref{eq:pair} and into~\eqref{eq:U-equation-derivative} gives
\begin{equation}\label{eq:linearAtZero}
\begin{aligned}
&cU_*\alpha+bU_*'\beta+aU_*''\delta=\frac5{18},\\
&c'U_*\alpha+(b'+c)U_*'\beta
 +\bigl(aU_*'''+(a'+b)U_*''\bigr)\delta=0.
\end{aligned}
\end{equation}
In~\eqref{eq:linearAtZero} and the remainder of this proof,
all functions and their derivatives are evaluated at the
first zero $q_0$.
Subtract $q\eta$ times the second equation in~\eqref{eq:linearAtZero}
from $(q\eta)'$ times the first equation in~\eqref{eq:linearAtZero}.
Using the definitions~\eqref{eq:K-definitions} for the coefficients
of $\alpha$ and $\delta$, we obtain
\begin{equation}\label{eq:weighted-difference}
\left[K_1+\frac5{18}(q\eta)'\right]\alpha
 +U_*'\bigl[b(q\eta)'-(b'+c)q\eta\bigr]\beta-K_2\delta
=\frac5{18}(q\eta)'.
\end{equation}
To simplify the coefficient of $\beta$, we differentiate the
second equation in~\eqref{eq:pair} and then evaluate at $q_0$:
\begin{equation}\label{eq:Ustar-equation-derivative}
aU_*'''+(a'+b)U_*''+(b'+c)U_*'+c'U_*=(1-p)(q\eta)'.
\end{equation}
Subtract $q\eta$ times equation~\eqref{eq:Ustar-equation-derivative}
from $(q\eta)'$ times the second equation in~\eqref{eq:pair}.
The result is
\[
\begin{aligned}
&U_*\bigl[c(q\eta)'-c'q\eta\bigr]
 +U_*'\bigl[b(q\eta)'-(b'+c)q\eta\bigr]\\
&\quad+U_*''\bigl[a(q\eta)'-(a'+b)q\eta\bigr]
 -aq\eta U_*'''=\frac5{18}(q\eta)',
\end{aligned}
\]
which implies, by~\eqref{eq:K-definitions}, that
\begin{equation}\label{eq:beta-coefficient}
U_*'\bigl[b(q\eta)'-(b'+c)q\eta\bigr]=K_2-K_1.
\end{equation}
Substituting~\eqref{eq:beta-coefficient}
into~\eqref{eq:weighted-difference} and collecting terms gives
\begin{equation}\label{eq:contradiction}
-\frac5{18}(q\eta)'(1-\alpha)+K_1(\alpha-\beta)+K_2(\beta-\delta)=0.
\end{equation}
By~\eqref{eq:ordering} and~\eqref{eq:wrange},
$1-\alpha$, $\alpha-\beta$, $\beta-\delta$ and $(q\eta)'$ are strictly positive.
Lemma~\ref{lem:signs} makes both $K_1$ and $K_2$ strictly negative.
Every term on the left of~\eqref{eq:contradiction} is therefore
strictly negative, a contradiction. Thus no first zero exists.
Together with \eqref{eq:Wcenter}, this proves the proposition.
\end{proof}

\section{Proofs of the main results}\label{sec:conclusion}

In this section we complete the proofs of the main results.

\subsection{The Chebyshev property}\label{subsec:chebyshev-proof}

We first prove Theorem~\ref{thm:Wronskian} by returning to the
variable $s$, and then apply Lemma~\ref{lem:ECT} to establish
the Chebyshev property of $(1,s,w,v)$.

\begin{proof}[Proof of Theorem~\ref{thm:Wronskian}]
For every $\kappa>1$ and $1<s<\kappa$, the corresponding
variables in~\eqref{eq:coordinates} satisfy $0<p<1$ and $0<q<1$.
By~\eqref{eq:conversion}, \eqref{eq:Wq} and
Proposition~\ref{prop:global},
\[
W(1,s,w,v)=\frac{(1-p)^5}{p^4}\mathscr W>0.\qedhere
\]
\end{proof}

\begin{proposition}\label{prop:ECT}
For every $\kappa>1$, the ordered tuple $(1,s,w(s),v(s))$
is an ECT-system on $(1,\kappa)$.
\end{proposition}
\begin{proof}
Its first two Wronskians equal one.
By~\eqref{eq:conversion} and Lemma~\ref{lem:ratio},
its third Wronskian is negative, and its fourth Wronskian
is positive by Theorem~\ref{thm:Wronskian}.
Thus all four successive Wronskians are nonzero on $(1,\kappa)$.
Lemma~\ref{lem:ECT} therefore implies that $(1,s,w,v)$ is an
ECT-system on this interval.
\end{proof}

\begin{corollary}\label{lem:zeros}
For every fixed $\kappa>1$, every linear combination of
$1$, $s$, $w(s)$ and $v(s)$ that is not identically zero
has at most three zeros on $(1,\kappa)$, counted with multiplicity.
\end{corollary}
\begin{proof}
Proposition~\ref{prop:ECT} shows that $(1,s,w,v)$
is an ECT-system on $(1,\kappa)$. By
Definition~\ref{def:Chebyshev} with $n=k=4$, every linear
combination that is not identically zero has at most
$k-1=3$ zeros, counted with multiplicity.
\end{proof}

\subsection{Zeros of the generating function}\label{subsec:generating-proof}

We first record the center expansion needed for the endpoint
version of Rolle's theorem.

\begin{lemma}\label{lem:I-center}
For every fixed $\kappa>1$, the function $I(h)$
in~\eqref{eq:I} extends analytically to $h=h_c$ and, as
$h\to h_c^+$, satisfies
\begin{equation}\label{eq:I-center}
I(h)=J_1(\kappa)
\left(-\frac23\mu_1+\mu_2+\mu_3-2(\kappa-1)\mu_4\right)
(h-h_c)+O\bigl((h-h_c)^2\bigr).
\end{equation}
In particular, $I(h_c)=0$.
\end{lemma}
\begin{proof}
By~\cite[Lemma~8 and its proof]{Zhao}, the functions
$I'_{i,j}(h(s))$ associated with the area integrals
in~\eqref{eq:I} are analytic at $s=\kappa$ and have the common
value $J_1(\kappa)$ there.
The proof of Lemma~8 in~\cite{Zhao} also gives
$I_{i,j}(h_c)=0$.
Using~\eqref{eq:s-variable} and integrating from $h_c$,
we obtain analytic extensions of these integrals satisfying
\[
I_{i,j}(h)=J_1(\kappa)(h-h_c)+O\bigl((h-h_c)^2\bigr)
\]
for the indices occurring in~\eqref{eq:I}.
Substituting these expansions into~\eqref{eq:I} and using
$h_c=-2/3$ yields~\eqref{eq:I-center}.
\end{proof}

We now combine Corollary~\ref{lem:zeros} with the center
condition in Lemma~\ref{lem:I-center} to complete the
estimate for the number of zeros of $I(h)$.

\begin{proof}[Proof of Theorem~\ref{thm:generating}]
By Lemma~\ref{lem:I-center}, we have $I(h_c)=0$.
If $\overline G(h)\equiv0$, identity~\eqref{eq:RolleI} makes
$I(h)/h$ constant. Its value at $h_c$ is zero, so $I(h)\equiv0$.
Thus $I(h)\not\equiv0$ implies $\overline G(h)\not\equiv0$.

By~\eqref{eq:G-ratios}, $G(s)/J_1(s)$ is a linear
combination of $1$, $s$, $w(s)$ and $v(s)$.
It is not identically zero because $\overline G(h)\not\equiv0$.
Corollary~\ref{lem:zeros} therefore bounds its number of zeros
by three, counted with multiplicity.
Multiplication by $J_1(s)>0$ and the change $h=h(s)$
in~\eqref{eq:s-variable} preserve these zeros and their
multiplicities, so $\overline G(h)$ also has at most three zeros.

Because $h\ne0$ on $[h_c,h_s)$, $I(h)$ and $I(h)/h$ have the same
interior zeros and multiplicities. Applying
Lemma~\ref{lem:endpoint} to $I(h)/h$ and using~\eqref{eq:RolleI}
shows that $I(h)$ has at most as many zeros on $(h_c,h_s)$ as
$\overline G(h)$, counted with multiplicity. Hence $I(h)$ has at
most three zeros there.
\end{proof}

\subsection{Cyclicity}

In this subsection, we prove the main results of this paper.

\begin{proof}[Proof of Theorem~\ref{thm:main}]
For a generic center, let $\Pi$ denote its period annulus and
use $\kappa>1$ as in Section~\ref{sec:reduction}.
By the reduction to the generating function~\eqref{eq:I}
described in Section~\ref{sec:reduction},
Theorem~\ref{thm:generating} implies
$\operatorname{Cycl}(\Pi)\leq3$. The cyclicity is at least three by~\cite[Theorem~1]{Zhao}.
Thus $\operatorname{Cycl}(\Pi)=3$.

For the non-generic center $b=2$, $c=0$, the cyclicity is
three by~\cite[Theorem~5]{Iliev1998}. For $b=-2$, $c=0$,
the same conclusion follows from~\cite[Theorem~1.2]{Iliev1997}.

The above three cases cover all centers in $Q_4$, completing the proof.
\end{proof}

Finally, we apply Proposition~\ref{prop:prescribed-zeros}
to the four-dimensional space of generating functions to
prove Theorem~\ref{thm:interpolation}.

\begin{proof}[Proof of Theorem~\ref{thm:interpolation}]
For the given generic center, let $h_1<h_2<h_3$ be the
first-integral values of the three prescribed periodic orbits.
By~\cite[Corollary~9(i)]{Zhao}, the generating functions
in~\eqref{eq:I} form a four-dimensional real vector space
of analytic functions on $(h_c,h_s)$.
By Theorem~\ref{thm:generating}, this space satisfies the
hypotheses of Proposition~\ref{prop:prescribed-zeros}, which
gives a function $I(h)$ whose zeros are precisely $h_1,h_2,h_3$,
all simple.
The coefficient realization in~\cite[Theorems~1 and~2(iii)]{Iliev1998}
and the transformations in~\cite{GI} realize this generating
function by an analytic quadratic perturbation.
For all sufficiently small nonzero parameter values, the
three simple zeros yield three hyperbolic limit cycles
converging to the prescribed periodic orbits.
\end{proof}

\section*{Acknowledgments}
This work is supported by the National Natural Science Foundation of China
(No. 12371183).

\end{document}